\documentclass[12pt,letterpaper, english]{article} 
\usepackage[activeacute,english]{babel}
\usepackage[]{inputenc}
\usepackage{moreverb}
\usepackage{natbib}
\usepackage{amsmath}
\usepackage{graphicx,color} 
\usepackage{amsfonts} 
\usepackage{amssymb}
\usepackage{amsthm}
\usepackage{amsxtra}
\usepackage{amsgen}
\usepackage{amscd}
\usepackage{tensor}
\usepackage{multirow}
\usepackage{latexsym}
 \usepackage{tabulary} 
\usepackage [colorlinks, bookmarks open, bookmarks numbered,citecolor= blue, urlcolor= blue]{hyperref} 
\usepackage{float}
\newcommand {\R}{\mathbb{R}}

\newcommand{\N}{\mathbb{N}}
\newcommand{\Z}{\mathbb{Z}}

\newtheorem{lemma}{\textbf{Lemma}}
\newtheorem{theorem}{\textbf{Theorem}}

\newtheorem{defi}{Definition}[section]

\newtheorem{afir}{Affirmation}[section] 
\date{} 

\begin{document}

           \title{A new iterative three-point method for solving systems of nonlinear equations.}
           \maketitle  
            
\begin{center}
\author{\textbf{ Carlos E. Cadenas R.$^{1,2}$, Yorman J. Mendoza N.$^{3}$}}\\           
{\scriptsize
 
 $^{1}$ Department of Mathematics, Experimental Faculty of Science and Technology, University of Carabobo, Venezuela.\\
 $^{2}$ Multidisciplinary Center for Visualization and Scientific Computing, University of Carabobo, Venezuela .\\
 $^{3}$ Department of Physics and Mathematics of the University National Experimental Francisco de Miranda, Venezuela .\\
}
\end{center}	
\begin{abstract}\par
A three-point iterative method for solving scalar non-linear equations was selected and then adapted to solve systems of non-linear equations. Subsequently, by applying Taylor's theorem to functions of $\R^{n}$ in $\R^{n}$, it is shown that the new method also has a sixth order of convergence. It is confirmed that the theoretical order of convergence coincides with the computational order of convergence by the numerical solution of two problems. Finally, its computational efficiency is calculated and subsequently compared with that of other three-point methods of fifth and sixth order convergence that also solve systems of non-linear equations.\\

\end{abstract}

\footnote{For email: Carlos E. Cadenas R, ccadenas@uc.edu.ve },\footnote{ccadenas45@gmail.com}
\footnote{ Yorman J. Mendoza N, myorman602@gmail.com, yjmsigma@gmail.com }.

\section{Introduction}
\par
The construction and analysis of iterative fixed-point methods for solving equations and systems of nonlinear equations is one of the most attractive topics in the theory of numerical analysis with wide applications in science, engineering, economics and finance. \cite{Sharma2016}. A system of non-linear equations is an equation of the form $\textbf{F}(\pmb{x})=\textbf{0}$ where $\textbf{F}:D\subseteq\R^{n}\longrightarrow\R^{n}$  is a nonlinear vector function $\textbf{F}(\pmb{x})=\left(f_{1}(\pmb{x}),f_{2}(\pmb{x}),\cdots,f_{n}(\pmb{x})\right)^{T}$, with $\pmb{x}=\left(x_{1},x_{2},\cdots,x_{n}\right)^{T}\in D$. Solving such a system consists of finding a vector $\pmb{\alpha}$ such that $\pmb{\alpha}=\left(\alpha _{1},\alpha_{2},\cdots,\alpha _{n}\right)^{T}$ and $\textbf{F}(\pmb{\alpha})=\textbf{0}$. The solution vector $\pmb{\alpha}$ of $\textbf{F}(\pmb{x})=\textbf{0}$ is obtained as a fixed point of some iteration function $\Phi:\R^{n}\longrightarrow\R^{n}$ sufficiently continuous, corresponding to the iteration equation:
\begin{equation}\label{Equapuntfijo}
\pmb{x}^{\left(k+1\right)} =\Phi\left(\pmb{x}^{\left(k\right)} \right),k=0,1,2,\cdots
\end{equation}
\par

Newton's method is one of the most commonly used single-point iterative techniques for solving systems of nonlinear equations \cite{Ostrowski1960}, \cite{Ortega1970},\cite{Kelly2003} and is given by
\begin{equation}\label{metodonewton} 
  \pmb{x}^{(k+1)}=\pmb{x}^{(k)}-\textbf{F}^{\,'}(\pmb{x}^{(k)})^{-1}\textbf{F}(\pmb{x}^{(k)}),     
\end{equation} 
where $\textbf{F}'(\textbf{x}^{(k)})$ is the first derivative of Fréchet and $\textbf{F}^{\,'}(\textbf{x}^{(k)})^{-1}$ is inverse in the $k$- iteration, having this quadratic convergence order and a small computational cost, since for each iteration it uses only one functional evaluation and one matrix inversion.
In recent years, a good number of publications relating to the solution of systems of nonlinear equations have been oriented towards the construction of methods and families of methods that improve the order of convergence of Newton's method, since there are applications in which large systems of nonlinear equations arise, where the computational speed with which a solution is obtained plays a fundamental role, being this closely related to the order of convergence. This has allowed the rise of multipoint methods to solve systems of nonlinear equations, because they overcome the limitations that one-point methods present with respect to the order of convergence and computational efficiency, \cite{Vassileva2011}. Among these we can mention \cite{Homeier2004}, \cite{Frontinini2004}, \cite{Cordero2006}, \cite{Darvishi2007-1}, \cite{Cordero2}, \cite{Noor2009}, \cite{Cordero3}, \cite{Miquel2011}, \cite{Cordero4}, \cite{cordero5}, \cite{Sharma2014}, \cite{Sharma20142}, \cite{Xiao2015}, \cite{Sharma3}, \cite{Esmaeili2015}, \cite{Sharma2016}, \cite{Narang2016},\cite{GowLorWink2024} among others.\par

The most direct way to build an iterative method for solving systems of nonlinear equations is to adapt methods that solve a nonlinear equation (one variable) to the vector case (several variables), but this implies greater difficulty in obtaining new iterative methods for solving systems of nonlinear equations (see \cite{Vassileva2011}), since there are inherent aspects of the constructed method that are directly inherited from the method that has been extended from the scalar case to the vector case.\\\par
In this sense, in the present article, a sixth-order three-point method of convergence is obtained for solving systems of nonlinear equations by extending to the vector case the three-point method for solving scalar nonlinear equations constructed by Parhi and Gupta \cite{ParhiGutpa2008}, whose iterative scheme is the following:

\begin{equation}\label{MetPG} 
 \begin{array}{rcl}
  y_n &=& x_n-\frac{f(x_n)}{f'(x_n)},\\\\ 
  z_n &=& x_n-\frac{2f(x_n)}{f'(x_n)+f'(y_n)},\\\\
  x_{n+1} &=& z_n-\frac{f(z_n)}{f'(x_n)}\frac{f'(x_n)+f'(y_n)}{3f'(y_n)-f'(x_n)},
  	\end{array}
\end{equation}

which uses two functional evaluations and two first derivatives for each iteration, and whose order of convergence is six. The rest of the article is organized as follows. Section 2 presents some definitions related to convergence. Section 3 presents the generalization to the vector case of the iterative scheme (\ref{MetPG}) and its respective convergence analysis. Section 4 confirms that the theoretical order of convergence coincides with the computational order of convergence by solving two numerical examples. Section 5 shows the computational efficiency of the constructed method and compares it with the efficiency of two fifth and sixth order three-point methods for systems of nonlinear equations given in \cite{Cordero4} and \cite{Cordero3}, respectively. In section 6 the conclusions are given.

\section{Definitions relating to convergence}
The concepts related to the order of convergence, approximation error and error equation are presented below:\\

\begin{defi}
Let $\left\{\textbf{x}^{(k)}\right\}_{k\geq 0}$ a succession in $\R^{n}$ that converges to $\pmb{\alpha}\in\R^{n}$. Then convergence is called:
 \begin{itemize}
	\item[a)] 
	\textbf{Linear}, if  there exists  $0<M<1$, $k_{0}\in\N$ such that :\\\\
	$\left\|\textbf{x}^{(k+1)}-\pmb{\alpha}\right\|\leq M\left\|\textbf{x}^{(k)}-\pmb{\alpha}\right\|$, $\forall k\geq k_{0}$.
	\item[b)]
	\textbf{Quadratic}, if there exists $M>0$, $k_{0}\in\N $ such that:\\\\
	$\left\|\textbf{x}^{(k+1)}-\pmb{\alpha}\right\|\leq M\left\|\textbf{x}^{(k)}-\pmb{\alpha}\right\|^{2}$, $\forall k\geq k_{0}$.
	
	\item[c)] \textbf{Of order $\textbf{p}$}, if there exists  $M>0$, $k_{0}\in \N $ such that:\\\\  
	$\left\|\textbf{x}^{(k+1)}-\pmb{\alpha}\right\|\leq M\left\|\textbf{x}^{(k)}-\pmb{\alpha}\right\|^{p}$, $\forall k\geq k_{0}$.
  \end{itemize}
\end{defi}.

\begin{defi}
Let $\left\{\textbf{x}^{(k)}\right\}_{k\geq 0}$ be a sequence in $\R^{n}$ that converges to $\pmb{\alpha}$, through an iterative process described by the equation(\ref{Equapuntfijo}). The expression given by $\textbf{e}_k=\textbf{x}^{(k)}-\pmb{\alpha}$ is called the \textbf{approximation error in the $k$-\'esima iteraci\'on} and the relationship:

\begin{equation}\label{error}
	\textbf{e}_{k+1}=L\left(\textbf{e}_k\right)^{p}+O\left(\left(\textbf{e}_k\right)^{p+1}\right)
\end{equation}
is defined as the \textbf{error equation}. Where $p$ is the order of convergence, $L$ a function $p$-linear bounded of $\overbrace{\R^{n}\times\R^{n}\times\cdots\times\R^{n}}^{\textbf{p}-times}$ in $\R^{n}$.
\end{defi}
\section{Development of the method}
 Generalizing to $\R^{n}$ the iterative scheme described by (\ref{MetPG}), We obtain a three-point iterative method for solving systems of nonlinear equations as follows:

\begin{equation}\label{PGvectorial} 
 \begin{array}{rcl}

  \textbf{y}^{(k)}&=&\textbf{x}^{(k)}-\textbf{F}^{\,'}(\textbf{x}^{(k)})^{-1}\textbf{F}(\textbf{x}^{(k)}),\\\\
	
  \textbf{z}^{(k)}&=&\textbf{x}^{(k)} -2\left[\textbf{F}^{\,'}(\textbf{x}^{(k)})+\textbf{F}^{\,'}(\textbf{y}^{(k)})\right]^{-1}\textbf{F}(\textbf{x}^{(k)}),\,k=0,1,2,\cdots \\\\
  \textbf{x}^{(k+1)}&=&\textbf{z}^{(k)} -\left[3\textbf{F}^{\,'}(\textbf{y}^{(k)})-\textbf{F}^{\,'}(\textbf{x}^{(k)})\right]^{-1}\left[\textbf{F}^{\,'}(\textbf{x}^{(k)})+\textbf{F}^{\,'}(\textbf{y}^{(k)})\right]\textbf{F}^{\,'}(\textbf{x}^{(k)})^{-1}\textbf{F}(\textbf{z}^{(k)}),
	\end{array}
\end{equation} 
the starting point $\textbf{x}^{(0)}$ is given, $\textbf{F}'(\textbf{x}^{(k)})$ is the first derivative of Fréchet and $\textbf{F}^{\,'}(\textbf{x}^{(k)})^{-1}$ is inverse in the $k$-\'esima iteration. To obtain the theoretical order of convergence of the three-point iterative scheme (\ref{PGvectorial}), we need the following result for functions of $\R^{n}$ in $\R^{n}$ ( see \cite{Ortega1970}).
 
\begin{lemma}
If $\textbf{F}:D\subseteq\R^{n}\longrightarrow\R^{n}$ is a function $p$-times Fréchet differentiable in the convex set $D\subseteq\R^{n}$ then for $\textbf{x},\,\textbf{y}\in D$ we have to:
\begin{equation}
\textbf{F}(\textbf{y})=\sum^{p-1}_{k=0} \frac{1}{k!}\textbf{F}^{\,\,(k)}(\textbf{x})(\textbf{y-x})^{k}+\textbf{R}_{p}.
\end{equation}

where:
\begin{equation}\nonumber
	\textbf{R}_{p}=\int^{1}_{0}\frac{(1-t)^{p-1}}{(p-1)!}\textbf{F}^{(p)}(\textbf{x}+t(\textbf{y-x}))(\textbf{y-x})^{p}dt
\end{equation}
 $$(\textbf{y-x})^{p}=\overbrace{\left((\textbf{y-x}),\,(\textbf{y-x}),\,\ldots,\,(\textbf{y-x})\right)}^{p-veces }$$, and 
\begin{equation}\nonumber
\left\|\textbf{R}_{p}\right\|\leq\frac{1}{p!}\sup_{0 \leq t \leq 1}\left\{\left\|\textbf{F}^{(p)}(\textbf{x}+t(\textbf{y-x}))\right\|\right\}\left\|\textbf{y-x}\right\|^{p}.
\end{equation}
\end{lemma}
The following theorem states that the order of convergence of the constructed iterative scheme is six.
  \begin{theorem}
\noindent Sea $\textbf{F}:D\subseteq\R^{n}\longrightarrow\R^{m}$ sufficiently Fréchet differentiable in $D$ and $\pmb{\alpha}\in\,D$ is a solution to the system of nonlinear equations $\textbf{F}(\textbf{x})=\textbf{0}$. Suppose that $\textbf{F}^{\,'}$ is continuous and non-singular on an open ball $B(\pmb{\alpha},r)\subset D$. So, the succession $\left\{\textbf{x}^{(k)}\right\}_{k\geq 0}$ obtained using the iterative scheme (\ref{PGvectorial}) converges to $\pmb{\alpha}$ with order of convergence six and its error equation is given by:

$$\textbf{e}_{k+1}=\left(9\textbf{A}_3\textbf{A}_2+{\frac{4}{3}}\textbf{A}^{3}_2-{\frac{34}{3}}\textbf{A}_2\textbf{A}_3-{\frac{3}{4}}\textbf{A}^{3}_2\textbf{A}_3-{\frac{2}{3}}\textbf{A}_4-6\textbf{A}_3\right)\left(\textbf{A}_3+{\frac{3}{4}}\textbf{A}^{2}_2\textbf{A}_3\right)\textbf{e}^{6}_k+O(\textbf{e}^{7}_k)$$.

\end{theorem}

\textbf{Proof}. First, we develop the Taylor series for $\textbf{F}$ and $\textbf{F}'$around $\pmb{\alpha}$ in $\textbf{x}^{(k)}\in B(\pmb{\alpha},r)\subset D$, So:\\

\begin{equation}\label{Fx}
\textbf{F}(\textbf{x}^{(k)})=\textbf{F}'(\pmb{\alpha})\left[\textbf{e}_k+\textbf{A}_2\textbf{e}^{2}_k+\textbf{A}_3\textbf{e}^{3}_k+\textbf{A}_4\textbf{e}^{4}_k+\textbf{A}_5\textbf{e}^{5}_k+\textbf{A}_6\textbf{e}^{6}_k+\textbf{A}_7\textbf{e}^{7}_k\right]+O(\textbf{e}^{8}_k).
\end{equation}\par
Con $\textbf{A}_{j}=\frac{1}{j!}\textbf{F}'(\pmb{\alpha})^{-1}\textbf{F}^{(j)}(\alpha)\in \textit{L}(\overbrace{(\R^{n}\times\R^{n}\times\ldots\times\R^{n})}^{j-veces},\R^{n})$, $j=2,3,\cdots$, $\textbf{F}(\pmb{\alpha})=\textbf{0}$. similarly:\\
\begin{equation}\label{F'x}
\textbf{F}'(\textbf{x}^{(k)})=\textbf{F}'(\pmb{\alpha})\left[\textbf{I}+2\textbf{A}_2\textbf{e}_k+3\textbf{A}_3\textbf{e}^{2}_k+4\textbf{A}_4\textbf{e}^{3}_k+5\textbf{A}_5\textbf{e}^{4}_k+6\textbf{A}_6\textbf{e}^{5}_k+7\textbf{A}_7\textbf{e}^{6}_k\right]+O(\textbf{e}^{7}_k).\\
\end{equation}\par
with $\textbf{I}$ the identity matrix of order $n\times n$. Given that $\textbf{F}'(\textbf{x}^{(k)})^{-1}\cdot\textbf{F}'(\textbf{x}^{(k)})=\textbf{I}$ and $$\textbf{F}'(\textbf{x}^{(k)})^{-1}=\left[\textbf{I}+2\textbf{A}_2\textbf{e}_k+3\textbf{A}_3\textbf{e}^{2}_k+4\textbf{A}_4\textbf{e}^{3}_k+5\textbf{A}_5\textbf{e}^{4}_k+6\textbf{A}_6\textbf{e}^{5}_k+7\textbf{A}_7\textbf{e}^{6}_k\right]^{-1}\textbf{F}'(\pmb{\alpha})^{-1}+O(\textbf{e}^{7}_k),$$ we calculate $\textbf{F}'(\textbf{x}^{(k)})^{-1}$ in the form:\\

\begin{equation}\label{inversaF'x}
\textbf{F}'(\textbf{x}^{(k)})^{-1}=\left[\textbf{X}_1+\textbf{X}_2\textbf{e}_k+\textbf{X}_3\textbf{e}^{2}_k+\textbf{X}_4\textbf{e}^{3}_k+\textbf{X}_5\textbf{e}^{4}_k+\textbf{X}_6\textbf{e}^{5}_k+\textbf{X}_7\textbf{e}^{6}_k\right]\textbf{F}'(\pmb{\alpha})^{-1}+O(\textbf{e}^{7}_k),
\end{equation}
where, $\textbf{X}_1=\textbf{I}$ y $\textbf{X}_s=-\sum^{s}_{i=2}i\textbf{X}_{s-i+1}\textbf{A}_i$, $s=2,3,\cdots$. Then performing the multiplications of the equations (\ref{inversaF'x}) and (\ref{Fx}), grouping terms, we have that:\\
$$\textbf{F}'(\textbf{x}^{(k)})^{-1}\textbf{F}(\textbf{x}^{(k)})=\textbf{e}_k-\textbf{A}_2\textbf{e}^{2}_k+\textbf{M}_1\textbf{e}^{3}_k+\textbf{M}_2\textbf{e}^{4}_k+\textbf{M}_3\textbf{e}^{5}_k+\textbf{M}_4\textbf{e}^{6}_k+\textbf{M}_5\textbf{e}^{7}_k+O(\textbf{e}^{8}_k).$$
whit $\textbf{M}_s=\textbf{A}_{s+2}+\sum^{s}_{j=1}\textbf{X}_{j+1}\textbf{A}_{s-(j-2)}+\textbf{X}_{s+2}$, $s=1,2,\cdots.$ So:
\begin{equation}
\begin{array}{rcl}\nonumber	
	\widetilde{\textbf{e}}_k&=&\textbf{y}^{(k)}-\pmb{\alpha}=\textbf{x}^{(k)}-\pmb{\alpha}-\textbf{F}'(\textbf{x}^{(k)})^{-1}\textbf{F}(\textbf{x}^{(k)})\\\\
    &=&\textbf{e}_k-\left[\textbf{e}_k-\textbf{A}_2\textbf{e}^{2}_k+\textbf{M}_1\textbf{e}^{3}_k+\textbf{M}_2\textbf{e}^{4}_k+\textbf{M}_3\textbf{e}^{5}_k+\textbf{M}_4\textbf{e}^{6}_k+\textbf{M}_5\textbf{e}^{7}_k+ O(\textbf{e}^{8}_k)\right].\\
\end{array} 
\end{equation}\par 

As a consequence:\\
\begin{equation}\label{Error1}
\widetilde{\textbf{e}}_k=\textbf{A}_2\textbf{e}^{2}_k-\textbf{M}_1\textbf{e}^{3}_k-\textbf{M}_2\textbf{e}^{4}_k-\textbf{M}_3\textbf{e}^{5}_k-\textbf{M}_4\textbf{e}^{6}_k-\textbf{M}_5\textbf{e}^{7}_k+O(\textbf{e}^{8}_k).
\end{equation}\par
Developing the Taylor series for $\textbf{F}$ and $\textbf{F}'$ around of $\pmb{\alpha}$ in the point $\textbf{y}^{(k)}\in B(\pmb{\alpha},r)\subset D$ we have:
$$\textbf{F}(\textbf{y}^{(k)})=\textbf{F}'(\pmb{\alpha})\left[\widetilde{\textbf{e}}_k+\textbf{A}_2\widetilde{\textbf{e}}^{2}_k+A_3\widetilde{\textbf{e}}^{3}_k\right]+O(\widetilde{\textbf{e}}^{4}), $$ and $$\textbf{F}'(\textbf{y}^{(k)})=\textbf{F}'(\pmb{\alpha})\left[\textbf{I}+2\textbf{A}_2\widetilde{\textbf{e}}_k+3\textbf{A}_3\widetilde{\textbf{e}}^{2}_k+4A_4\widetilde{\textbf{e}}^{3}_k\right]+O(\widetilde{\textbf{e}}^{4}).$$\par

Substituting the equation (\ref{Error1}) in the previous equations respectively, developing the necessary notable products, grouping terms and applying the properties of the Landau $O$ function, we have that:\\ 
\begin{equation}\nonumber
\textbf{F}(\textbf{y}^{(k)})=\textbf{F}'(\pmb{\alpha})\left[\textbf{Q}_2\textbf{e}^{2}_k+2\textbf{Q}_3\textbf{e}^{3}_k+\textbf{Q}_4\textbf{e}^{4}_k+\textbf{Q}_5\textbf{e}^{5}_k+5\textbf{Q}_6\textbf{e}^{6}_k+\textbf{Q}_7\textbf{e}^{7}_k\right]+O(\textbf{e}^{8}_k).
\end{equation}
where:
\begin{equation}
	\begin{array}{rcl}\nonumber
	\textbf{Q}_2&=&\textbf{A}_2\\\\
	\textbf{Q}_3&=&-\textbf{M}_1\\\\
	\textbf{Q}_4&=&-\textbf{M}_2+\textbf{A}^{3}_2\\\\
	\textbf{Q}_5&=&-\textbf{M}_3-2\textbf{A}^{2}_2\textbf{M}_1\\\\
	\textbf{Q}_6&=&-\textbf{M}_4-2\textbf{A}^{2}_2\textbf{M}_2+\textbf{A}_2\textbf{M}^{2}_1+\textbf{A}_3\textbf{A}^{3}_2\\\\
	\textbf{Q}_7&=&-\textbf{M}_5-2\textbf{A}_2\textbf{M}_3+2\textbf{A}_2\textbf{M}_1\textbf{M}_2-3\textbf{A}_3\textbf{A}^{2}_2\textbf{M}_1,\\
	\end{array}
\end{equation}\par

Also:
\begin{equation}\label{F'y}
\textbf{F}'(\textbf{y}^{(k)})=\textbf{F}'(\pmb{\alpha})\left[\textbf{I}+\textbf{T}_2\textbf{e}^{2}_k+\textbf{T}_3\textbf{e}^{3}_k+\textbf{T}_4\textbf{e}^{4}_k+\textbf{T}_5\textbf{e}^{5}_k+\textbf{T}_6\textbf{e}^{6}_k+\textbf{T}_7\textbf{e}^{7}_k\right]+O(\textbf{e}^{8}_k).
\end{equation}\par
where:
\begin{equation}\nonumber
	\begin{array}{rcl}
	\textbf{T}_2&=&2\textbf{A}^{2}_2\\\\
	\textbf{T}_3&=&-2\textbf{A}_2\textbf{M}_1\\\\
	\textbf{T}_4&=&-2\textbf{A}_2\textbf{M}_2+3\textbf{A}_3\textbf{A}^{2}_2\\\\
	\textbf{T}_5&=&-2\textbf{A}_3\textbf{M}_3-6\textbf{A}_3\textbf{A}_2\textbf{M}_1\\\\
	\textbf{T}_6&=&-2\textbf{A}_2\textbf{M}_4+3\textbf{A}_3\textbf{M}_1-6\textbf{A}_3\textbf{A}_2\textbf{M}_2+4\textbf{A}_4+\textbf{A}^{3}_2\\\\
	\textbf{T}_7&=&-2\textbf{A}_2\textbf{M}_5-6\textbf{A}_3\textbf{A}_2\textbf{M}_3+6\textbf{A}_3\textbf{M}_1\textbf{M}_2-12\textbf{A}_4\textbf{A}^{2}_2\textbf{M}_1.\\
	\end{array}
\end{equation}\par
From the equations (\ref{F'x}) y (\ref{F'y}): 
\begin{subequations}
\begin{equation}\label{F'ymasF'x}
\textbf{F}'(\textbf{y}^{(k)})+\textbf{F}'(\textbf{x}^{(k)})=\textbf{F}'(\pmb{\alpha})\left[2\textbf{I}+2\textbf{A}_2\textbf{e}_k+\textbf{L}_2\textbf{e}^{2}_k+\textbf{L}_3\textbf{e}^{3}_k+\textbf{L}_4  \textbf{e}^{4}_k+\textbf{L}_5\textbf{e}^{5}_k+\textbf{L}_6\textbf{e}^{6}_k\right]+O(\textbf{e}^{7}_k).
\end{equation}
\begin{equation}\label{tresF'ymenosF'x}
3\textbf{F}'(\textbf{y}^{(k)})-\textbf{F}'(\textbf{x}^{(k)})=\textbf{F}'(\pmb{\alpha})\left[2\textbf{I}-2\textbf{A}_2\textbf{e}_k+\textbf{N}_2\textbf{e}^{2}_k+\textbf{N}_3\textbf{e}^{3}_k+\textbf{N}_4\textbf{e}^{4}_k+\textbf{N}_5\textbf{e}^{5}_k+\textbf{N}_6\textbf{e}^{6}_k\right]+O(\textbf{e}^{7}_k).
\end{equation} 
\end{subequations}
Where $\textbf{L}_s=\textbf{T}_s+(s+1)\textbf{A}_{s+1}$ and $\textbf{N}_s=3\textbf{T}_s-(s+1)\textbf{A}_{s+1}$, $s=2,3,\cdots$.\\\par
knowing that:
$$\left[\textbf{F}'(\textbf{y}^{(k)})+\textbf{F}'(\textbf{x}^{(k)})\right]^{-1}\left[\textbf{F}'(\textbf{y}^{(k)})+\textbf{F}'(\textbf{x}^{(k)})\right]=\textbf{I}$$ and 
 
$$\left[\textbf{F}'(\textbf{y}^{(k)})+\textbf{F}'(\textbf{x}^{(k)})\right]^{-1}=\left[2\textbf{I}+2\textbf{A}_2\textbf{e}_k+\textbf{L}_2\textbf{e}^{2}_k+\textbf{L}_3\textbf{e}^{3}_k+\textbf{L}_4  \textbf{e}^{4}_k+\textbf{L}_5\textbf{e}^{5}_k+\textbf{L}_6\textbf{e}^{6}_k\right]^{-1}\textbf{F}'(\pmb{\alpha})^{-1}+O(\textbf{e}^{7}_k),$$ 
We calculate $\left[\textbf{F}'(\textbf{y}^{(k)})+\textbf{F}'(\textbf{x}^{(k)})\right]^{-1}$ in the form:

\begin{equation}\nonumber
\left[\textbf{F}'(\textbf{y}^{(k)})+\textbf{F}'(\textbf{x}^{(k)})\right]^{-1}=\left[\textbf{Y}_1+\textbf{Y}_2\textbf{e}_k+\textbf{Y}_3\textbf{e}^{2}_k+\textbf{Y}_4\textbf{e}^{3}_k+\textbf{Y}_5  \textbf{e}^{4}_k+\textbf{Y}_6\textbf{e}^{5}_k+\textbf{Y}_7\textbf{e}^{6}_k\right]\textbf{F}'(\pmb{\alpha})^{-1}+O(\textbf{e}^{7}_k),
\end{equation}
in which,$\textbf{Y}_1=\frac{1}{2}\textbf{I}$, $\textbf{Y}_2=-\frac{1}{2}\textbf{A}_2$, and $\textbf{Y}_s=-\frac{1}{2}\left(2\textbf{Y}_{s-1}\textbf{A}_2+\sum^{s-2}_{i=1}\textbf{Y}_i\textbf{L}_{s-i}\right)$, $s=3,4,\cdots$.\\\par
Consequently:
\begin{equation}\label{invF'ymasF'x}
\left[\textbf{F}'(\textbf{y}^{(k)})+\textbf{F}'(\textbf{x}^{(k)})\right]^{-1}=\left[\frac{1}{2}\textbf{I}-\frac{1}{2}\textbf{A}_2\textbf{e}_k+\textbf{Y}_3\textbf{e}^{2}_k+\textbf{Y}_4\textbf{e}^{3}_k+\textbf{Y}_5  \textbf{e}^{4}_k+\textbf{Y}_6\textbf{e}^{5}_k+\textbf{Y}_7\textbf{e}^{6}_k\right]\textbf{F}'(\pmb{\alpha})^{-1}+O(\textbf{e}^{7}_k).
\end{equation}\par 
Similarly, since:$$\left[3\textbf{F}'(\textbf{y}^{(k)})-\textbf{F}'(\textbf{x}^{(k)})\right]^{-1}\left[3\textbf{F}'(\textbf{y}^{(k)})-\textbf{F}'(\textbf{x}^{(k)})\right]=\textbf{I},$$ 
we search $\left[3\textbf{F}'(\textbf{y}^{(k)})-\textbf{F}'(\textbf{x}^{(k)})\right]^{-1}$ in the form:

$$\left[3\textbf{F}'(\textbf{y}^{(k)})-\textbf{F}'(\textbf{x}^{(k)})\right]^{-1}=\left[\textbf{W}_1+\textbf{W}_2\textbf{e}_k+\textbf{W}_3\textbf{e}^{2}_k+\textbf{W}_4\textbf{e}^{3}_k+\textbf{W}_5  \textbf{e}^{4}_k+\textbf{W}_6\textbf{e}^{5}_k+\textbf{W}_7\textbf{e}^{6}_k\right]\textbf{F}'(\pmb{\alpha})^{-1}+O(\textbf{e}^{7}_k),$$
con $\textbf{W}_1=\frac{1}{2}\textbf{I},\,\textbf{W}_2=\frac{1}{2}\textbf{A}_2$, and $\textbf{W}_s=\textbf{W}_{s-1}\textbf{A}_2-\sum^{s-2}_{i=1}\textbf{W}_i\textbf{N}_{s-i}$, $s=3,4,\cdots$. Consequently:
\begin{equation}\label{invtresF'ymenosF'x}
\left[3\textbf{F}'(\textbf{y}^{(k)})-\textbf{F}'(\textbf{x}^{(k)})\right]^{-1}=\left[\frac{1}{2}\textbf{I}+\frac{1}{2}\textbf{A}_2\textbf{e}_k+\textbf{W}_3\textbf{e}^{2}_k+\textbf{W}_4\textbf{e}^{3}_k+\textbf{W}_5  \textbf{e}^{4}_k+\textbf{W}_6\textbf{e}^{5}_k+\textbf{W}_7\textbf{e}^{6}_k\right]\textbf{F}'(\pmb{\alpha})^{-1}+O(\textbf{e}^{7}_k).
\end{equation}\par
Multiplying the equations (\ref{invF'ymasF'x}) and (\ref{Fx}), so: 

\begin{equation}\nonumber
2\left[\textbf{F}'(\textbf{y}^{(k)})+\textbf{F}'(\textbf{x}^{(k)})\right]^{-1}\textbf{F}(\textbf{x}^{(k)})=\textbf{e}_k-\textbf{C}\textbf{e}^{3}_k-\textbf{R}_4\textbf{e}^{4}_k-\textbf{R}_5\textbf{e}^{5}_k-\textbf{R}_6\textbf{e}^{6}_k+O(\textbf{e}^{7}_k).
\end{equation}
with $\textbf{C}=-(\textbf{A}_3-\textbf{A}^2_2\textbf{Y}_3)$, $\textbf{R}_s=\textbf{A}_s-\textbf{A}_2\textbf{A}_{s-1}+2\left(\sum^{s-1}_{i=3}\textbf{Y}_i \textbf{A}_{s-(i-1)}+\textbf{Y}_s\right),\,s=4,5,\cdots.$ So:
\begin{equation}
	\begin{array}{rcl}\nonumber
	\widehat{\textbf{e}}_k&=&\textbf{z}^{(k)}-\pmb{\alpha}=\textbf{x}^{(k)}-\pmb{\alpha}-2\left[\textbf{F}'(\textbf{y}^{(k)})+\textbf{F}'(\textbf{x}^{(k)})\right]^{-1}\textbf{F}(\textbf{x}^{(k)})\\\\\nonumber
	&=&\textbf{e}_k-\left[\textbf{e}_k-\textbf{C}\textbf{e}^{3}_k-\textbf{R}_4\textbf{e}^{4}_k-\textbf{R}_5\textbf{e}^{5}_k-\textbf{R}_6\textbf{e}^{6}_k+O(\textbf{e}^{7}_k)\right].\\
		\end{array} 
\end{equation}\par
therefore:
\begin{equation}\label{Error2}
 \widehat{\textbf{e}}_k=\textbf{C}\textbf{e}^{3}_k+\textbf{R}_4\textbf{e}^{4}_k+\textbf{R}_5\textbf{e}^{5}_k+\textbf{R}_6\textbf{e}^{6}_k+O(\textbf{e}^{7}_k).
\end{equation}\par
Developing the Taylor series for $\textbf{F}$ around of $\pmb{\alpha}$ in the point $\textbf{z}^{(k)}\in B(\pmb{\alpha},r)\subset D$ then:
$$\textbf{F}(\textbf{z}^{(k)})=\textbf{F}'(\pmb{\alpha})\left[\widehat{\textbf{e}}_k+\textbf{A}_2\widehat{\textbf{e}}^{2}_k+\textbf{A}_3\widehat{\textbf{e}}^{3}_k\right]+O(\widehat{\textbf{e}}^{4}_k).$$

Substituting the equation (\ref{Error2}) in the previous equation, developing the necessary products, grouping terms and applying the properties of the Landau $O$ function we have that:\\
\begin{equation}\label{Fz}
	\textbf{F}(\textbf{z}^{(k)})=\textbf{F}'(\pmb{\alpha})\left[\textbf{C}\textbf{e}^{3}_k+\textbf{R}_4\textbf{e}^{4}_k+\textbf{R}_5\textbf{e}^{5}_k+\pmb{\lambda}\textbf{e}^{6}_k+\pmb{\beta}\textbf{e}^{6}_k\right]+O(\textbf{e}^{8}_k),\,
\end{equation}
where:
\begin{equation}
	\begin{array}{rcl}
	\pmb{\lambda}&=&\textbf{R}_6+\textbf{A}_2\textbf{C}^{2}\\\\\nonumber
	 \pmb{ \beta}&=&\textbf{R}_7+2\textbf{A}_2\textbf{C}\,\textbf{R}_4+2\textbf{A}_2\textbf{C}\,\textbf{R}_5. \\
	\end{array} 
\end{equation}\par
Multiplying the equations (\ref{inversaF'x}) and (\ref{Fz}):
\begin{equation}\label{EQa}
	\textbf{F}'(\textbf{x}^{(k)})^{-1}\textbf{F}(\textbf{z}^{(k)})=\textbf{C}\textbf{e}^{3}_k+\textbf{U}_4\textbf{e}^{4}_k+\textbf{U}_5\textbf{e}^{5}_k+\textbf{U}_6\textbf{e}^{6}_k+\textbf{U}_7\textbf{e}^{7}_k+O(\textbf{e}^{8}_k),
\end{equation}
With:
\begin{equation}\nonumber
	\begin{array}{rcl}
	\textbf{U}_4&=&\textbf{R}_4-2\textbf{A}_2\textbf{C}\\\\\nonumber
	\textbf{U}_5 &=&\textbf{R}_5-2\textbf{A}_2\textbf{R}_4+\textbf{X}_3\textbf{C}\\\\\
	\textbf{U}_6 &=&\pmb{\lambda}-2\textbf{A}_2\textbf{R}_5+\textbf{X}_3\textbf{R}_4+\textbf{X}_4\textbf{C}\\\\\
		\textbf{U}_7 &=&\pmb{\beta}-2\textbf{A}_2\pmb{\lambda}+\textbf{X}_3\textbf{R}_5+\textbf{X}_5\textbf{C}.
	\end{array} 
\end{equation}\par
  
Multiplying the equations (\ref{invtresF'ymenosF'x}) and (\ref{F'ymasF'x}, so:
 \begin{equation}\label{EQb}
\left[3\textbf{F}'(\textbf{y}^{(k)})-\textbf{F}'(\textbf{x}^{(k)})\right]^{-1}\left[\textbf{F}'(\textbf{y}^{(k)})+\textbf{F}'(\textbf{x}^{(k)})\right]=\textbf{I}+2\textbf{A}_2\textbf{e}_k+\textbf{E}_2\textbf{e}^{2}_k+\textbf{E}_3\textbf{e}^{3}_k+\textbf{E}_4\textbf{e}^{4}_k+\textbf{E}_5\textbf{e}^{5}_k+\textbf{E}_6\textbf{e}^{6}_k+O(\textbf{e}^{7}_k),
\end{equation}
where:
\begin{equation}
	\begin{array}{rcl}
	\textbf{E}_2&=&\frac{1}{2}\textbf{L}_2+\textbf{A}^{2}_2+2\textbf{W}_3\\\\\nonumber
	\textbf{E}_3&=&\frac{1}{2}\textbf{L}_3+\frac{1}{2}\textbf{A}_2\textbf{L}_2+2\textbf{W}_3\textbf{A}_2+2\textbf{W}_4\\\\\
	\textbf{E}_s &=&\frac{1}{2}\textbf{L}_{s}+\frac{1}{2}\textbf{A}_2\textbf{L}_{s-1}+2\textbf{W}_{s}\textbf{A}_2+\sum^{s-1}_{i=3}\textbf{W}_{i}\textbf{L}_{s+1-i}+2\textbf{W}_s,\,s=4,6,\cdots.
	\end{array} 
\end{equation}\par
Multiplying the equations(\ref{EQb}) and (\ref{EQa})then:
 \begin{equation}
\left[3\textbf{F}'(\textbf{y}^{(k)})-\textbf{F}'(\textbf{x}^{(k)})\right]^{-1}\left[\textbf{F}'(\textbf{y}^{(k)})+\textbf{F}'(\textbf{x}^{(k)})\right]\textbf{F}'(\textbf{x}^{(k)})^{-1}\textbf{F}(\textbf{z}^{(k)})=\textbf{C}\textbf{e}^{3}_k+\Phi_4\textbf{e}^{4}_k+\Phi_5\textbf{e}^{5}_k+\Phi_6\textbf{e}^{6}_k+\Phi_7\textbf{e}^{7}_k+O(\textbf{e}^{8}_k),
\end{equation}

with:
\begin{equation}
	\begin{array}{rcl}
	\Phi_4&=&\textbf{U}_4+2\textbf{A}_2\textbf{C}\\\\\nonumber
	\Phi_5 &=&\textbf{U}_5+2\textbf{A}_2\textbf{U}_4+\textbf{E}_2\textbf{C}\\\\\
	\Phi_6 &=&\textbf{U}_6+2\textbf{A}_2\textbf{U}_5+\textbf{E}_2\textbf{U}_4+\textbf{E}_3\textbf{C}\\\\\
	\Phi_7 &=&\textbf{U}_7+2\textbf{A}_2\textbf{U}_6+\textbf{E}_2\textbf{U}_5+\textbf{E}_3\textbf{U}_4+\textbf{E}_4\textbf{C},
	\end{array} 
\end{equation}
so:
\begin{equation}\nonumber
	\begin{array}{rcl}
	\textbf{e}_{k+1}&=&\textbf{x}^{(k+1)}-\pmb{\alpha}=\textbf{z}^{(k)}-\pmb{\alpha}-\left[3\textbf{F}'(\textbf{y}^{(k)})-\textbf{F}'(\textbf{x}^{(k)})\right]^{-1}\left[\textbf{F}'(\textbf{y}^{(k)})+\textbf{F}'(\textbf{x}^{(k)})\right]\textbf{F}'(\textbf{x}^{(k)})^{-1}\textbf{F}(\textbf{z}^{(k)})\\\\
	&=& \widehat{\textbf{e}}_k-\left[3\textbf{F}'(\textbf{y}^{(k)})-\textbf{F}'(\textbf{x}^{(k)})\right]^{-1}\left[\textbf{F}'(\textbf{y}^{(k)})+\textbf{F}'(\textbf{x}^{(k)})\right]\textbf{F}'(\textbf{x}^{(k)})^{-1}\textbf{F}(\textbf{z}^{(k)})\\\\
	&=&(\textbf{R}_4-\Phi_4)\textbf{e}^{4}_k+(\textbf{R}_5-\Phi_5)\textbf{e}^{5}_k+(\textbf{R}_6-\Phi_6)\textbf{e}^{6}_k+O(\textbf{e}^{7}_k).
		\end{array} 
\end{equation}\par
Since $\Phi_4=\textbf{U}_4+2\textbf{A}_2\textbf{C}$ y $\textbf{U}_4=\textbf{R}_4-2\textbf{A}_2\textbf{C}$, we have $\Phi_4=\textbf{R}_4$ and $\textbf{R}_4-\Phi_4=\textbf{0}$. Similarly: 
$$\textbf{R}_5-\Phi_5=4\textbf{A}^{2}_2\textbf{C}-(\textbf{X}_3+\textbf{E}_2)\textbf{C},\,\textbf{X}_3=-3\textbf{A}_3+4\textbf{A}^{2}_2,\,\textbf{E}_2=3\textbf{A}_3.$$\par

So, $\textbf{R}_5-\Phi_5=4\textbf{A}^{2}_2\textbf{C}-(-3\textbf{A}_3+4\textbf{A}^{2}_2+3\textbf{A}_3)\textbf{C}=\textbf{0}$. Finally $\textbf{e}_{k+1}=(\textbf{R}_6-\Phi_6)\textbf{e}^{6}_k+O(\textbf{e}^{7}_k).$\par
Replacing the values of $\textbf{R}_6$, $\Phi_6$ and simplifying we have that:
$$\textbf{e}_{k+1}=\left(9\textbf{A}_3\textbf{A}_2+{\frac{4}{3}}\textbf{A}^{3}_2-{\frac{34}{3}}\textbf{A}_2\textbf{A}_3-{\frac{3}{4}}\textbf{A}^{3}_2\textbf{A}_3-{\frac{2}{3}}\textbf{A}_4-6\textbf{A}_3\right)\left(\textbf{A}_3+{\frac{3}{4}}\textbf{A}^{2}_2\textbf{A}_3\right)\textbf{e}^{6}_k+O(\textbf{e}^{7}_k)$$.

\section{Numerical tests}
Next, the constructed method is then used to solve two systems of nonlinear equations $\textbf{F}(\pmb{x})=\pmb{0}$. Numerical calculations have been performed using the software \textbf{Matlab R2019a} (\textbf{MAtrix LABoratory}) considering 20 significant digits.\\\par The stopping criterion used is $\left\|\pmb{x}^{(k+1)}-\pmb{x}^{(k)}\right\|<tol$, where $tol=10^{-12}$. In each problem we analyze the number of iterations necessary for the method to converge to the solution of the system of nonlinear equations and the theoretical order of convergence is confirmed with the computational order of convergence $p$ approximated by the equation (see \cite{Cordero3}): 
$$p\approx \frac{\ln\left(\left\|\pmb{x}^{(k+1)}-\pmb{x}^{(k)}\right\|/\left\|\pmb{x}^{(k)}-\pmb{x}^{(k-1)}\right\|\right)}{\,\,\,\ln\left(\left\|\pmb{x}^{(k)}-\pmb{x}^{(k-1)}\right\|/\left\|\pmb{x}^{(k-1)}-\pmb{x}^{(k-2)}\right\|\right)},$$.\\\par
The value considered for the order of computational convergence is the maximum of the $p$ values obtained.\\
\textbf{Problem 1}. Let us consider the system of two equations
$$\left\{
 \begin{array}{rcl}
x_{1}^{3}x_{2}^{3}-1&=&0,\\\
& & \\\
x_{1}-1&=&0.
\end{array} \right.$$ \par 

Taking the initial approximation $\pmb{x}^{(0)}=\left(2,2\right)^{T}$ the solution is $\pmb{\alpha}=\left(1,1\right)^{T}$ as shown in the table 3.1.\\\par
 
\begin{table}[h!]
\centering
\scalebox{0.8}[0.8]{\begin{tabular}{lllll}
           & & & &   \\\hline
$k$ &\,\,\,\,\,\,\,\,\,\,\,\,\,\,\,\,\,\,\,\,\,\,\,\,\,\,\,\,\,\,$x^{(k)}_1$&\,\,\,\,\,\,\,\,\,\,\,\,\,\,\,\,\,\,\,\,\,\,\,\,\,\,\,\,\,\,\,$x^{(k)}_2$&\,\,\,\,\,\,\,\,\,\,\,\,\,\,\,\,\,\,$\left\|\pmb{x}^{(k+1)}-\pmb{x}^{(k)}\right\|$&\,\,\,\,\,\,\,$p$\\\hline 
 0&2.000000000000000000000000     &2.000000000000000000000000&   &\\
 1&1.000000000000000000000000     &2.276866652619243600000000& 1.276866652619243800000000&\\
 2&1.000000000000000000000000     &1.041198047519996700000000&0.041198047519996617000000 &\\
 3&1.000000000000000000000000     &1.000000000869789100000000&0.000000000869789072852836 &5.9950\\ 
 4&1.000000000000000000000000     &1.000000000000000000000000&0.000000000000000000000000 & \\ \hline
\end{tabular}} 
\caption{Problem results 1.}\label{tabla1}
\end{table} 
\textbf{Problem 2}. Now, a problem of border values ( see \cite{Ortega1970}):
$$y''=\frac{1}{2}y^{3}+3y'-\frac{3}{2-x}+\frac{1}{2},\,y(0)=0,\,y(1)=1,$$
partitioning the interval $\left[0,1\right]$, such that $0=x_0<x_1<\cdots<x_n=1$, where $x_{i+1}=x_i+h$, $h=\frac{1}{n}$, then:
$$0=y(0)=y(x_0),y_1=y(x_1),y_2=y(x_2),\cdots,y_n=y(x_n)=1.$$\par
Discretizing the problem using the numerical approximation formulas for the first and second derivatives,

$$y'_j=\frac{y_{j+1}-y_{j-1}}{2h},\,y''_j=\frac{y_{j-1}-2y_j+y_{j+1}}{h^{2}},\,j=1,2,\cdots,n-1,$$

we have :
\begin{eqnarray*}
	\frac{y_{j-1}-2y_j+y_{j+1}}{h^{2}}&=&\frac{1}{2}y^{3}_j+3\left(\frac{y_{j+1}-y_{j-1}}{2h}\right)-\frac{3}{2-x_j}+\frac{1}{2}\\\\
	y_{j-1}-2y_j+y_{j+1}&=&\frac{\,\,h^2}{2}y^{3}_j+\frac{3h}{2}\left(y_{j+1}-y_{j-1}\right)-\frac{3h^{2}}{2-x_j}+\frac{\,\,h^2}{2}.
\end{eqnarray*}\par 
This is how a non-linear system is constructed $\left(n-1\right)$ equations with $\left(n-1\right)$ unknowns given by:
\begin{equation}
 \left(1+\frac{3}{2}h\right)y_{j-1}+\left(1-\frac{3}{2}h\right)y_{j+1}-2y_j-\frac{h^{2}}{2}y^{3}_j+\frac{\,3h^{2}}{2-x_j}-\frac{h^{2}}{2}=0,\,j=1,2,\cdots,n-1.\end{equation}\par   
Particularly for $n=7$, $y_0=y(0)=0$, $y_7=y(x_7)=1$:
\begin{eqnarray*}
	f_1(y_1,y_2,\cdots,y_6)&=&-2y_1-\frac{\,\,h^2}{2}y^{3}_1+(1-\frac{3}{2}h)y_2+\frac{\,3h^{2}}{2-x_1}-\frac{\,\,h^2}{2}\\\\
	f_j(y_1,y_2,\cdots,y_6)&=&\left(1+\frac{3}{2}h\right)y_{j-1}+\left(1-\frac{3}{2}h\right)y_{j+1}-2y_j-\frac{\,\,h^{2}}{2}y^{3}_j+\frac{\,3h^{2}}{2-x_j}-\frac{h^{2}}{2},\,j=2,3,4,5\\\\
	f_6(y_1,y_2,\cdots,y_6)&=&\left(1+\frac{3}{2}h\right)y_{5}-2y_6-\frac{\,\,h^{2}}{2}y^{3}_6+\frac{\,3h^{2}}{2-x_6}+\frac{\,\,h^{2}}{2}+\left(1-\frac{3}{2}h\right).
\end{eqnarray*}\par 
This non-linear system is reduced to solving: $$\textbf{F}(y_1,y_2,\cdots,y_6)=\left(f_1(y_1,y_2,\cdots,y_6),f_2(y_1,y_2,\cdots,y_6),\cdots,f_6(y_1,y_2,\cdots,y_6)\right)^{T}=\pmb{0}.$$
 For $\pmb{x}^{(0)}=\left( 7{.}25,7{.}25,7{.}25,7{.}25,7{.}25,7{.}25\right)^{T}$
                                    , we have $\pmb{\alpha}\approx\begin{pmatrix}
                                     \,\,\,\,0{.}0083494505929842810\\
                                     \,\,\,\,0{.}0077178744477447688\\
                                     \,\,\,\,0{.}0257257242738628680\\
                                     -0{.}0169564948809037690\\
                                     \,\,\,\,0{.}1244784293587575000\\
																		 -0{.}2954656773600666300\\
                               \end {pmatrix}$, as shown in the tables (\ref{tabla2}) y (\ref{tabla3}).

\begin{table}[h!]
\centering
\scalebox{0.7}[0.7]{\begin{tabular}{lllll}
           & & & &   \\\hline
$k$ &\,\,\,\,\,\,\,\,\,\,\,\,\,\,\,\,\,\,\,\,\,\,\,\,\,\,\,\,\,\,$x^{(k)}_1$&\,\,\,\,\,\,\,\,\,\,\,\,\,\,\,\,\,\,\,\,\,\,\,\,\,\,\,\,\,\,\,$x^{(k)}_2$&\,\,\,\,\,\,\,\,\,\,\,\,\,\,\,\,\,\,\,\,\,\,\,\,\,\,\,\,\,\,\,$x^{(k)}_3$&\,\,\,\,\,\,\,\,\,\,\,\,\,\,\,\,\,\,\,\,\,\,\,\,\,\,\,\,\,\,\,$x^{(k)}_4$\\\hline 
 0&\,\,7.250000000000000000000000     &\,\,7.250000000000000000000000&\,\,7.250000000000000000000000 &7.250000000000000000000000\\
 1&-0.017114396804569942000000    & -0.059580334143502722000000&  -0.055541505243312250000000& -0.158837854766547220000000\\
 2&\,\,0.008349450592892408400000     & \,\,0.007717874447977668400000&\,\, 0.025725724273191825000000 &-0.016956494878929907000000\\
 3&\,\,0.008349450592984281000000    & \,\,0.007717874447744768800000 &\,\, 0.025725724273862868000000& -0.016956494880903769000000\\ 
 4&\,\,0.008349450592984281000000     & \,\,0.007717874447744768800000& \,\,0.025725724273862868000000&-0.016956494880903769000000\\ \hline
\end{tabular}} 
\caption{Problem Results 2.}\label{tabla2}
\end{table}
 
\begin{table}[h!] 
\centering
\scalebox{0.8}[0.8]{\begin{tabular}{lllll}
           & & & &   \\\hline
$k$ &\,\,\,\,\,\,\,\,\,\,\,\,\,\,\,\,\,\,\,\,\,\,\,\,\,\,\,\,\,\,\,$x^{(k)}_5$&\,\,\,\,\,\,\,\,\,\,\,\,\,\,\,\,\,\,\,\,\,\,\,\,\,\,\,\,\,\,\,$x^{(k)}_6$&\,\,\,\,\,\,\,\,\,\,\,\,\,\,\,\,\,\,$\left\|\pmb{x}^{(k+1)}-\pmb{x}^{(k)}\right\|$&\,\,\,\,\,$p$\\\hline 
 0&7.250000000000000000000000&\,\,7.250000000000000000000000&&\\
 1&0.023407390755885604000000& -0.513568512472107660000000&18.083123849280749000000000 &\\
 2&0.124478429356049410000000& -0.295465677343173920000000&\,\,0.299493473744831430000000&\\
 3&0.124478429358757500000000& -0.295465677360066630000000&\,\,0.000000000017236776477730& \,\,5.7499\\ 
 4&0.124478429358757500000000& -0.295465677360066630000000&\,\,0.000000000000000000000000&\\ \hline
\end{tabular}} 
\caption{Continuation of results of the Problem 2.}\label{tabla3}
\end{table} 			

\section{Computational efficiency} 
To measure the computational efficiency of this iterative method, the Ostrowski efficiency index will be used \cite{Ostrowski1960}, which establishes that the computational efficiency of an iterative method is given by $E=p^{\frac{1}{C}}$, where $p$ is the order of convergence and $C$ the computational cost. According to \cite{Miquel2011}, for an iterative method that solves systems of $n$ nonlinear equations with $n$ unknowns the computational cost is given by:
$$
C(\mu_{0},\mu_{1},n)=P_{0}(n)\mu_{0}+P_{1}(n)\mu_{1}+P(n),
$$
where $P_{0}(n)$ the number of evaluations of scalar functions $\left(f_{1},f_{2},\cdots,f_{n}\right)$used to evaluate $\textbf{F}$, $P_{1}(n)$ the number of evaluations of scalar functions in $\textbf{F}^{\,'}$,that is to say $\frac{\partial f_{i}}{\partial x_{j}}$ with $1\leq i,j\leq n$, $P(n)$ the number of products and quotients needed per iteration and $\mu_{0}$, $\mu_{1}$ are the proportions between products and evaluations required to express the value of $C(\mu_{0},\mu_{1},n)$ in terms of products.\\\par       

For calculate $\textbf{F}$ in an iterative method, are evaluated $n$ scalar functions, for each new derivative $\textbf{F}^{\,'}$ are evaluated $n^{2}$ scalar functions, Additionally, the computational work to calculate the inverse of a matrix where products of the form appear is included $\pmb{A}^{-1}\cdot \pmb{b}$. Now well, instead of calculating the inverse, the linear system is solved $\pmb{A} \cdot \pmb{y}=\pmb{b}$ where are relized $n(n-1)(2n-1)/6$ products and $n(n-1)/2$ quotients in the decomposition $LU$, and $n(n-1)$ products and $n$ quotients in the resolutions of  two triangular linear systems, adicionally $n^{2}$ products in the  multiplication of a matrix with a vector or of a scalar by a matrix and $n$ products for the multiplication of a scalar by a vector. In addition to this it is assumed that a quotient is equivalent to $l$ products.\\\par
Now we denote by $\pmb{\phi}^{(6)}_{1}$ the iterative scheme which is our object of study, so the method efficiency index $\pmb{\phi}^{(p)}_{i}$ is  $E^{(p)}_{i}$ and  its computational cost $C^{(p)}_{i}$, where $p$ is the order of convergence of the $i$-th method. In this sense in $\pmb{\phi}^{(6)}_{1}$ for each iteration it is evaluated $\textbf{F}$ twice, therefore $P_{0}(n)=2n$, $\textbf{F}^{\,'}$ is evaluated twice, then $P_{1}(n)=2n^{2}$. Three inverse matrices are calculated, for what they do $3n(n-1)(2n-1)/6$ products and $3n(n-1)/2$ quotients in the decomposition $LU$ and $3n(n-1)$ products and $3n$ quotients in the resolution of six triangular linear systems (two for each inverse matrix), is added $6n^{2}$ products, from the five matrix-vector products and a scalar-matrix product existing in each iteration. lastly they are considered $n$ products originating from the multiplication of a scalar by a vector. In this sense, we have that:

	\begin{eqnarray*}
	P(n)&=&\frac{3n(n-1)(2n-1)}{6}+3n(n-1)+6n^{2}+n+l\left(\frac{3n(n-1)}{2}+3n\right)\\\\
	&=&\frac{n(n-1)(2n-1)}{2}+3n(n-1)+6n^{2}+n+l\frac{n}{2}\left(3(n-1)+6\right) \\\\
	&=&\frac{n}{2}\left(2n^{2}+15n-3+3l(n+1)\right).
\end{eqnarray*} \par

   Therefore:
\begin{equation}\label{CC1}
			C^{(6)}_{1}=2n\mu_{0}+2n^{2}\mu_{1}+\frac{n}{2}\left(2n^{2}+15n-3+3l(n+1)\right)\,and \,E^{(6)}_{1}=6^{1/C^{(6)}_{1}}.
\end{equation}\par
			
Now the computational efficiency of $\pmb{\phi}^{(6)}_{1}$ It is compared with that of two fifth and sixth order three-point methods given by Cordero et al in \cite{Cordero4} and \cite{Cordero3} which we denote by $\pmb{\phi}^{(5)}_{1}$ and $\pmb{\phi}^{(6)}_{2}$ respectively.\newline \par
These methods are given by:

\textbf{Method} $\pmb{\phi}^{(5)}_{1}$: 

\begin{equation}\nonumber 
 \begin{array}{rcl}

  \textbf{y}^{(k)}&=&\textbf{x}^{(k)}-\textbf{F}^{\,'}(\textbf{x}^{(k)})^{-1}\textbf{F}(\textbf{x}^{(k)}),\\\\
	
  \textbf{z}^{(k)}&=&\textbf{x}^{(k)} -2\left[\textbf{F}^{\,'}(\textbf{x}^{(k)})+\textbf{F}^{\,'}(\textbf{y}^{(k)})\right]^{-1}\textbf{F}(\textbf{x}^{(k)}),\\\\
  \textbf{x}^{(k+1)}&=&\textbf{z}^{(k)}-\textbf{F}^{\,'}(\textbf{y}^{(k)})^{-1}\textbf{F}(\textbf{z}^{(k)}).
	\end{array}
\end{equation}

\textbf{Method} $\pmb{\phi}^{(6)}_{2}$: 

\begin{equation}\nonumber 
 \begin{array}{rcl}

  \textbf{y}^{(k)}&=&\textbf{x}^{(k)}-\frac{2}{3}\textbf{F}^{\,'}(\textbf{x}^{(k)})^{-1}\textbf{F}(\textbf{x}^{(k)}),\\\\

  \textbf{z}^{(k)}&=&\textbf{x}^{(k)} -\frac{1}{2}\left[3\textbf{F}^{\,'}(\textbf{y}^{(k)})-\textbf{F}^{\,'}(\textbf{x}^{(k)})\right]^{-1}\left[3\textbf{F}^{\,'}(\textbf{y}^{(k)})+\textbf{F}^{\,'}(\textbf{x}^{(k)})\right]\textbf{F}^{\,'}(\textbf{x}^{(k)})^{-1}\textbf{F}(\textbf{x}^{(k)}),\\\\
  \textbf{x}^{(k+1)}&=&\textbf{z}^{(k)} -2\left[3\textbf{F}^{\,'}(\textbf{y}^{(k)})-\textbf{F}^{\,'}(\textbf{x}^{(k)})\right]^{-1}\textbf{F}(\textbf{z}^{(k)}).
	\end{array}
\end{equation}\par
Following a procedure similar to that described to obtain the computational cost and efficiency of the method $\pmb{\phi}^{(6)}_{1}$, the efficiency and computational cost of the methods are obtained $\pmb{\phi}^{(5)}_{1}$ y $\pmb{\phi}^{(6)}_{2}$ which are given by:
\begin{eqnarray}
C^{(5)}_{1}&=&2n\mu_{0}+2n^{2}\mu_{1}+\frac{n}{2}\left(2n^{2}+9n-3+3l(n+1)\right)\, and \,E^{(5)}_{1}=5^{1/C^{(5)}_{1}}.\\\label{CC2}
C^{(6)}_{2}&=&2n\mu_{0}+2n^{2}\mu_{1}+\frac{n}{3}\left(2n^{2}+18n+4+3l(n+1)\right)\, and \,E^{(6)}_{2}=6^{1/C^{(6)}_{2}}.\label{CC3}
\end{eqnarray}

\subsection{Comparison of efficiencies}
To compare the computational efficiency of two iterative methods $\pmb{\phi}^{(p)}_{i}$ and $\pmb{\phi}^{(q)}_{j}$ the quotient or ratio is considered:
$$R^{p,q}_{i,j}=\frac{\log E^{(p)}_{i}}{\log E^{(q)}_{j}}=\frac{C^{q}_{j}\log(p)}{C^{p}_{i}\log(q)}.$$\par
Clearly, if  $R^{p,q}_{i,j}>1$ then $\pmb{\phi}^{(p)}_{i}$ is more efficient that $\pmb{\phi}^{(q)}_{j}$, by if $R^{p,q}_{i,j}<1$ the opposite happens, that is to say $\pmb{\phi}^{(p)}_{i}$ is less efficient that  $\pmb{\phi}^{(q)}_{j}$. The boundary between two efficiencies is given by $R^{p,q}_{i,j}=1$ which would be the starting point in the comparisons. When comparing methods of different order of convergence, it is possible to write $\mu_0$ as a function of $\mu_1$, $n$ and  $l$, where $\left(\mu_0,\mu_1\right)\in \left(0,+\infty\right)\times \left(0,+\infty\right)$, $n\in\Z^{+}$ is such that $n\geq 2$ and $l\geq 1$.\\\par 

\textbf{Comparing} $\pmb{\phi}^{(6)}_{1}$\textbf{with} $\pmb{\phi}^{(5)}_{1}$:\\
 In this case we consider $s=\log(6)$ and $t=\log(5)$, of equality $R^{p,q}_{i,j}=1$ is expressed in $\mu_0$ as a function of $\mu_1$, $n$ and $l$, in the following way: 
\begin{equation}\label{BC} 
1=R^{6,5}_{1,1}=\frac{C^{(5)}_{1}\log(6)}{C^{(6)}_{1}\log(5)}\Leftrightarrow C^{(6)}_{1}\log(5)=C^{(5)}_{1}\log(6)
\end{equation}\par
\begin{center}
\begin{figure}[ht!]
  \centering
  \includegraphics[scale=0.5]{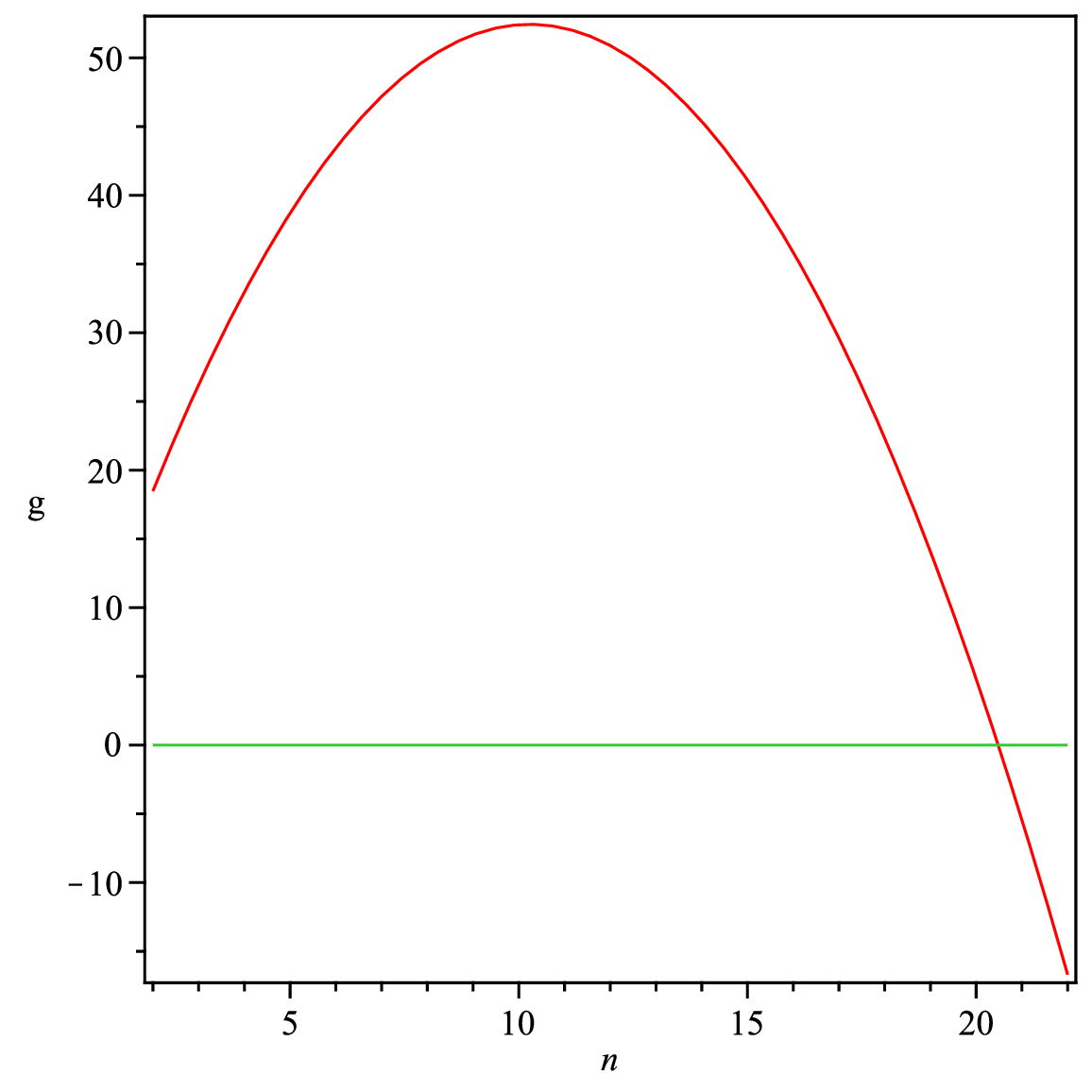}\newline \newline \newline \newline
	\caption{graph of $g$.}\label{eficiencia4}
\end{figure}
\end{center}

From the equations (\ref{CC1}), (\ref{CC2}) and (\ref{BC})one has:  
$$\mu_0=\frac{4n(t-s)\mu_1+2n^{2}(t-s)+3n(5t-3s)+3(s-t)+3n(t-s)l+3(t-s)l}{4(s-t)}=g(n,\mu_1,l),$$
considering $\mu_1>0$ and $l\geq1$ fixed, it has to be $g$ is always negative for $n\geq 21$, as shown in the figure \ref{eficiencia4}.

\begin{afir}
$\pmb{\phi}^{(6)}_{1}$ is more efficient than $\pmb{\phi}^{(5)}_{1}$ for $n\geq 21$.
\end{afir} 
In effect:\\

\begin{equation}\label{BC2}
1<R^{6,5}_{1,1}=\frac{C^{(5)}_{1}\log(6)}{C^{(6)}_{1}\log(5)}\Leftrightarrow C^{(6)}_{1}\log(5)<C^{(5)}_{1}\log(6),
\end{equation}
From the equations (\ref{CC1}), (\ref{CC2}) and (\ref{BC2}), have  $\mu_0>g(n,\mu_1,l)$, which is always true for $n\geq 21$, $\mu_1>0$ and $l\geq1$, since for those values $g$ is always negative and $\mu_0 \in \left(0,+\infty\right)$, for such reasons it is concluded that $\pmb{\phi}^{(6)}_{1}$ is more efficient than $\pmb{\phi}^{(5)}_{1}$ for $n\geq 21$.\\\par
 
\textbf{Comparing} $\pmb{\phi}^{(6)}_{1}$ \textbf{with} $\pmb{\phi}^{(6)}_{2}$:\\
In this case we compare the method obtained with another of the same order of convergence, To do this, it is sufficient to compare the computational cost of both, as shown: 
\begin{eqnarray*}
1=R^{6,6}_{1,2}=\frac{C^{(6)}_{2}\log(6)}{C^{(6)}_{1}\log(6)}& \Leftrightarrow & C^{(6)}_{1}=C^{(6)}_{2}\\\
& \Leftrightarrow &\frac{n}{2}\left(2n^{2}+15n-3+3l(n+1)\right)=\frac{n}{3}\left(2n^{2}+18n+4+3l(n+1)\right)\\\
&\Leftrightarrow&2n^{2}+9n+3l(n+1)-17=0.
\end{eqnarray*}\par
The equation $2n^{2}+9n+3l(n+1)-17=0$, gives us the border, however if $n\geq 2$ and $l\geq 1$ then $2n^{2}\geq 8$, $9n\geq 18$ and $3l(n+1)\geq 9$.\newline\par
\begin{center}

\begin{figure}[ht!]
  \centering
  \includegraphics[scale=0.8]{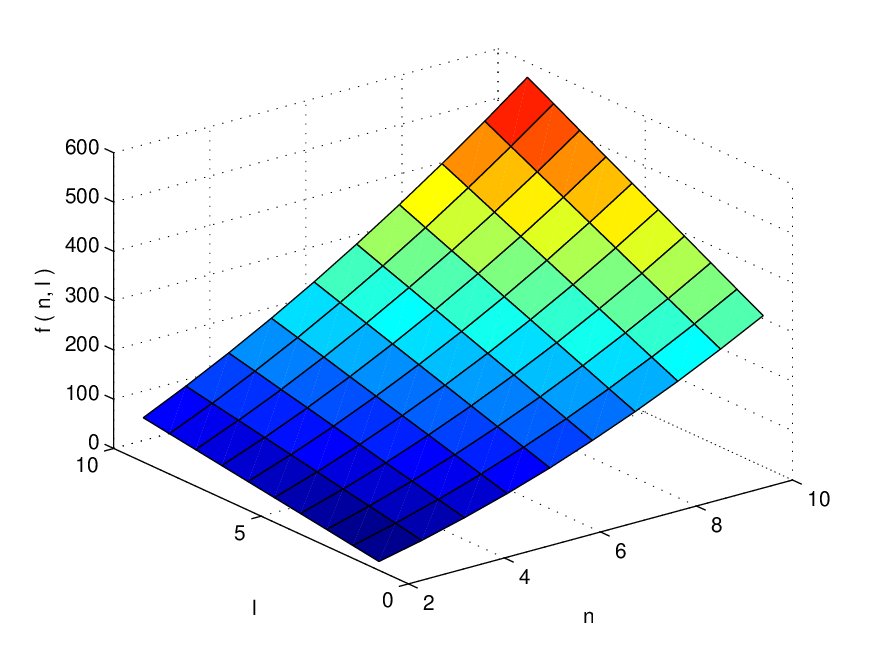}\newline \newline \newline  \newline  \newline \newline \newline \newline \newline \newline \newline \newline  \newline  \newline \newline\newline
	\caption{Graph of $f(n,l)=2n^{2}+9n+3l(n+1)-17$.}\label{eficiencia3}
\end{figure}
\end{center}

 In consecuently  $2n^{2}+9n+3l(n+1)\geq 8+18+9=35$ and  $2n^{2}+9n+3l(n+1)-17\geq 18>0$ as shown in the figure \ref{eficiencia3}.
\par
In this sense it is concluded that $ C^{(6)}_{1}>C^{(6)}_{2}$ for each $n\geq 2$ and $l\geq1$, consequently we have \newline $E^{(6)}_{1}=6^{1/C^{(6)}_{1}}<E^{(6)}_{2}=6^{1/C^{(6)}_{2}}$.

\newpage

\section{Conclusions}
In the present investigation, a three-point iterative method was constructed to solve systems of nonlinear equations, by extending to the vector case of a three-point iterative method to solve a nonlinear scalar equation. On the other hand, the respective analysis of it was given, based firstly on demonstrating analytically that the resulting order of convergence is six, secondly; its computational efficiency was calculated, proving to be more efficient than one of the two methods with which it was compared for the case of systems of nonlinear equations of size equal to or greater than 21 equations with 21 unknowns, it was observed that the number of matrix inversions negatively affects the computational efficiency of a method since it brings with it the largest set of computational operations. Furthermore, it was numerically verified that the computational order of convergence coincides with the analytical order of convergence. In this paper, the demonstration of convergence has been made from a local point of view. It should be noted that after repeated attempts to calculate the theoretical order of convergence of the method constructed with the help of the symbolic calculus software Maple and the impossibility of this to distinguish non-commutative multiplications of the matrix coefficients that arise in the Taylor expansion in functions of $\R^{n}$ in $\R^{n}$, the convergence calculations were made without the aid of any software, verifying through the exchange of communication via email with other researchers in the area, that the software \textbf{Mathematica} also presents the same drawback..\\\par
In the future the intention is to study the dynamics of the constructed method and to build new methods for systems of nonlinear equations by applying some of the techniques of numerical analysis that allow to reduce the number of inversions of Jacobian matrices and functional evaluations present in the iterative scheme, which increases the efficiency of the method and also reduces and facilitates the calculations to obtain the theoretical order of convergence, aspects that cannot be controlled when extending to the vector case existing methods that solve nonlinear scalar equations.
\cite{Vassileva2011}

\end{document}